\documentclass[twoside,11pt]{article}
\title{} \author{} \date{}
\usepackage{latexsym,amssymb,times,graphicx}
\input amssym.def
\newtheorem{te}{Theorem}[section]
\newtheorem{prop}[te]{Proposition}
\newtheorem{cor}[te]{Corollary}
\newtheorem{fac}[te]{Fact}

\newtheorem{rem}[te]{Remark}
\newtheorem{ex}[te]{Example}

\def\dok{\noindent{\bf Proof. }}
\def\kdok{\hfill $\Box$ \par \vspace*{2mm} }

\let\strokel\l

\def\a{\alpha}
\def\b{\beta}
\def\g{\gamma}
\def\d{\delta}

\def\o{\omega}
\def\k{\kappa}
\def\l{\lambda}
\def\r{\rho}
\def\s{\sigma}

\def\S{{\mathbb S}}
\def\T{{\mathbb T}}

\def\Q{{\mathbb Q}}
\def\B{{\mathbb B}}
\def\N{{\mathbb N}}
\def\X{{\mathbb X}}
\def\Y{{\mathbb Y}}
\def\Z{{\mathbb Z}}

\def\A{{\mathbb A}}

\def\BG{{\mathbb G}}
\def\BL{{\mathbb L}}

\def\CA{{\mathcal A}}
\def\CC{{\mathcal C}}

\def\L{{\mathcal L}}
\def\M{{\mathcal M}}
\def\CX{{\mathcal X}}
\def\CZ{{\mathcal Z}}
\def\CW{{\mathcal W}}
\def\c{{\mathfrak{c}}}
\def\la{\langle}
\def\ra{\rangle}

\def\id{\mathop{\mathrm{id}}\nolimits}

\def\otp{\mathop{\rm otp}\nolimits}
\def\PO{\mathop{\rm PO}\nolimits}
\def\LO{\mathop{\rm LO}\nolimits}
\def\Scatt{\mathop{\mbox{\rm Scatt}}\nolimits}
\def\SScatt{\mathop{\s\mbox{\rm -Scatt}}\nolimits}
\def\Card{\mathop{\rm Card}\nolimits}
\def\Ord{\mathop{\rm Ord}\nolimits}
\def\Mono{\mathop{\rm Mono}\nolimits}
\def\Cond{\mathop{\rm Cond}\nolimits}
\def\Emb{\mathop{\rm Emb}\nolimits}
\def\Iso{\mathop{\rm Iso}\nolimits}
\def\Aut{\mathop{\rm Aut}\nolimits}
\def\Sur{\mathop{\rm Sur}\nolimits}
\def\Sym{\mathop{\rm Sym}\nolimits}
\def\ar{\mathop{\rm ar}\nolimits}
\def\Mod{\mathop{\rm Mod}\nolimits}
\def\rfm{\mathop{\rm RFM}\nolimits}
\def\CR{\mathop{\mathcal{R}}\nolimits}
\def\wfr{\mathop{\rm Wfr}\nolimits}

\begin{document}
\thispagestyle{plain}
\begin{center}
           {\large \bf \uppercase{Reversibility of disconnected structures}}
\end{center}
\begin{center}
{\bf Milo\v s S.\ Kurili\'c\footnote{Department of Mathematics and Informatics, Faculty of Science, University of Novi Sad,
              Trg Dositeja Obradovi\'ca 4, 21000 Novi Sad, Serbia.
              email: milos@dmi.uns.ac.rs}
and Nenad Mora\v ca\footnote{Department of Mathematics and Informatics, Faculty of Science, University of Novi Sad,
              Trg Dositeja Obradovi\'ca 4, 21000 Novi Sad, Serbia.
              email:
              nenad.moraca@dmi.uns.ac.rs}}
\end{center}
\begin{abstract}
\noindent
A relational structure is called reversible iff every bijective endomorphism of that structure is an automorphism.
We give several equivalents of that property in the class of disconnected binary structures and some its subclasses.
For example, roughly speaking and denoting the set of integers by $\Z$, a structure having reversible components is reversible iff
its components can not be ``merged" by condensations (bijective homomorphisms)
and each $\Z$-sequence of condensations between different components must be, in fact, a sequence of isomorphisms.
We also give equivalents of  reversibility in some special classes of structures.
For example, we characterize CSB linear orders of a limit type
and show that a disjoint union of such linear orders is a reversible poset
iff the corresponding sequence of order types is finite-to-one.

\end{abstract}
\section{Introduction}\label{S1}
A relational structure $\X$ is said to be reversible iff every bijective endomorphism $f:\X\rightarrow \X$ is an automorphism
and the relevance of that property follows from the fact that
the class of reversible structures includes linear orders, Boolean lattices, well founded posets with finite levels \cite{Kuk,Kuk1},
tournaments, Henson graphs \cite{KuMoExtr}, and Henson digraphs \cite{KRet}.
In addition, reversibility is an invariant of some forms of bi-interpretability \cite{KRet},
extreme elements of $L_{\infty \o}$-definable classes of structures are reversible under some syntactical restrictions \cite{KuMoExtr},
and all structures first-order definable in linear orders by quantifier-free formulas without parameters (i.e., monomorphic or chainable structures)
are reversible \cite{KDef}.

In this article we investigate reversibility in the class of {\it binary structures}, that is models of the relational language
$L_b=\la R \ra$, where $\ar (R)=2$, and, moreover, we restrict our attention to the class of disconnected $L_b$-structures.
(If $\X =\la X,\r \ra$ is an $L_b$-structure and $\sim _\r$ the minimal equivalence relation on $X$
containing $\r $, then the corresponding equivalence classes are called the {\it connectivity components} of $\X$ and
$\X$ is said to be {\it disconnected} if it has more than one component, that is, if $\sim _\r \;\neq X^2$). The prototypical
disconnected structures are, of course, equivalence relations themselves; other prominent representatives of that class are
some countable ultrahomogeneous graphs and posets (see \cite{Lach,Sch}), non-rooted trees, etc.

By \cite{KuMoEqu}, a disconnected $L_b$-structure $\X =\bigcup _{i\in I}\X_i$
belonging to a wide class RFM\footnote{$\rfm$ is the class of
structures $\X =\bigcup _{i\in I}\X_i$ such that $\la \X_i :i\in I\ra$ is a
sequence of pairwise disjoint, connected and reversible $L_b$-structures
 which is {\it rich for monomorphisms}, which means that
for all $i,j\in I$ and each $A\in [X_j ]^{|X_i|}$ there is a monomorphism $g:\X _i \rightarrow \X _j$ such that $g[X_i]=A$.
} of structures (containing equivalence relations) is reversible iff
the sequence of cardinalities of its connectivity components, $\la |X_i|:i\in I\ra$, has the following property:
a sequence of non-zero cardinals $\la \k _i :i\in I\ra$ is defined to be
{\it reversible} iff
\begin{equation}\label{EQB032}\textstyle
\neg \exists f\in \Sur (I)\setminus \Sym (I)\;\;\forall j\in I \;\;\sum _{i\in f^{-1}[\{ j \}]}\k_i =\k _j ,
\end{equation}
where $\Sym (I)$ (resp.\ $\Sur (I)$) denotes the set of all bijections (resp.\ surjections) $f:I\rightarrow I$.
For the following 
characterization of such sequences we recall that a set $K\subset\N$ 
is called {\it independent}
iff $n\not\in \la K\setminus \{ n \}\ra$, for all $n\in K$,
where $\la K\setminus \{ n \}\ra$ is the subsemigroup of the semigroup  $\la \N , +\ra$
generated by $K\setminus \{ n \}$;  by $\gcd (K)$  we denote the greatest common divisor of the numbers from $K$. By \cite{KuMoEqu} we have
\begin{fac}\label{TB068}
A sequence of non-zero cardinals $\la \k _i \!:i\in I\ra$ is reversible if and only if
either it is a finite-to-one sequence, or it is a sequence of natural numbers, the set
$K=\{ m\in \N : \k _i =m ,\mbox{ for infinitely many } i\in I \}$
is non-empty and independent,
and $\gcd (K)$  divides at most finitely many elements
of the set $\{ \k _i :i\in I \}$.
\end{fac}
In Section \ref{S2} we give several equivalents of reversibility in the class of disconnected $L_b$-structures and here we describe one of them.
First, by Theorem \ref{TB066}, the components of a reversible structure must be reversible and, hence,
the following assumption (containing that restriction) will appear in several parts of the text:
\begin{itemize}
\item[($\ast$)] $\X _i =\la X _i ,\r _i\ra $, $i\in I$, are pairwise disjoint, connected and reversible $L_b$-structures, $\X = \bigcup _{i\in I}\X_i$
                and $\CX = \{ \X _i: i\in I \}$.
\end{itemize}
Second, generalizing (\ref{EQB032}) and writing $\X \preccurlyeq _c \Y$ iff there is a bijective homomorphism $g: \X \rightarrow \Y$,
we will say that a sequence of $L_b$-structures $\la \X _i :i\in I\ra$ is a
{\it reversible sequence of structures} iff
\begin{equation}\label{EQB040}\textstyle
\neg \exists f\in \Sur (I)\setminus \Sym (I)\;\;\forall j\in I \;\;\bigcup _{i\in f^{-1}[\{ j \}]}\X_i \preccurlyeq_c \X _j .
\end{equation}
Third, assuming ($\ast$), a mapping $i:\Z \rightarrow I$, usually denoted by $\la i_k :k\in \Z\ra$, will be called a
{\it $\Z$-sequence in $I$} iff it is an injection and
\begin{equation}\label{EQB041}\textstyle
\forall k,l\in \Z \;\;( k<l \Rightarrow \X_{i_k}\preccurlyeq_c\X_{i_l}),
\end{equation}
which is, by the transitivity
of the relation $\preccurlyeq_c$,
equivalent to the existence of a sequence of condensations  $g_k:\X_{i_{k}}\rightarrow\X_{i_{k+1}}$, for $k\in \Z$.
If, in addition, $\X_{i_k}\cong \X_{i_l}$, for all $k,l\in \Z$, the $\Z$-sequence $\la i_k :k\in \Z\ra$ will be called {\it trivial}.\footnote{We note that
$\la i_k :k\in \Z\ra$ is a $\Z$-sequence in $I$
iff $k\mapsto \X _{i_k}$ is a monomorphism from the linear order $\la \Z ,<\ra$ to the preorder $\la \CX,\preccurlyeq_c\ra$.
It is trivial iff all $\X _{i_k}$'s are in the same $\sim _c$-class (or, equivalently, $\cong$-class) of that preorder (see Fact \ref{TA011}).}

The announced equivalent of reversibility is given in the following theorem.
\begin{te}\label{TB065}
($\ast$) The structure $\X$ is reversible iff $\la \X _i :i\in I\ra$ is a reversible sequence of structures and
each $\Z$-sequence in $I$ is trivial.\footnote{In other words,
the union $\bigcup _{i\in I}\X_i$ is not reversible iff the components can be ``merged by condensations"
or there is a non-trivial $\Z$-sequence in $I$.}
\end{te}
In Section \ref{S2}, in addition,
we give equivalents of  reversibility in some special classes of structures (structures having finitely many components,
structures with finite components, disjoint unions of linear orders and tournaments).

In Sections \ref{S3} and \ref{S4} we give some sufficient conditions for reversibility of disconnected structures
and detect several classes of reversible partial orders. In Section \ref{S5} we apply these results to unions of disjoint linear orders.
In particular, we characterize CSB linear orders of a limit type
and show that a disjoint union of such linear orders, $\bigcup _{i\in I}\X_i$, is a reversible poset
iff $\la \otp (\X _i) :i\in I\ra$ is a finite-to-one sequence. Similarly, a disjoint union of $\s$-scattered linear orders, $\bigcup _{i\in I}\X_i$,
is reversible if $\la [\X _i]_\rightleftarrows :i\in I\ra$ is a finite-to-one sequence, where $[\X _i]_\rightleftarrows$ is the equimorphism
(bi-embedability) class of $\X _i$ in $\CX$.

In the rest of this section we introduce notation and recall basic facts which will be used in the paper.
\paragraph{Condensational order, equivalence and reversibility}
If $\X$ and $\Y$ are $L_b$-struc\-tures, $\Iso (\X ,\Y )$, $\Cond (\X ,\Y )$, $\Mono (\X ,\Y )$ and $\Emb (\X ,\Y )$
will denote the set of all isomorphisms, condensations (bijective homomorphisms), monomorphisms and embeddings from $\X$ to $\Y$ respectively.
$\Iso (\X ,\X )=\Aut (\X )$ is the set of automorphisms of $\X$, instead of $\Cond (\X ,\X )$ we will write $\Cond (\X )$ etc.

The {\it condensational preorder} $\preccurlyeq _c $ on the class of $L_b$-structures is defined by $\X \preccurlyeq _c \Y$ iff $\Cond (\X ,\Y )\neq\emptyset$,
the {\it condensational equivalence} is the equivalence relation defined on the same class by: $\X \sim _c \Y$ iff
$\X \preccurlyeq _c \Y$ and $\Y \preccurlyeq _c \X$,
and it determines the antisymmetric quotient of the condensational preorder, the {\it condensational order}, in the usual way.
We will write $\X \prec _c \Y$ iff $\X \preccurlyeq _c \Y$ and $\Y \not\sim _c \X$
(which is for reversible structures equivalent to $\X \preccurlyeq _c \Y$ and $\Y \not\cong \X$, see Fact \ref{TA011}).

Some simple properties equivalent to reversibility are listed in the following claim (which, in fact, holds for any relational language $L$, see \cite{KuMoVar}).
\begin{fac}\label{TB045}
For an $L_b$-structure $\X =\la X,\r\ra$ the following conditions are equivalent

(a) $\X $ is a reversible structure (that is, $\Cond (\X )=\Aut (\X )$),

(b) $\forall\s\subset X^2 \;\;(\s\varsubsetneq \r \Rightarrow \la X,\s\ra\ \not\cong \la X,\r\ra )$,

(c) $\forall\s\subset X^2 \;\;(\r\varsubsetneq \s  \Rightarrow \la X,\s\ra\ \not\cong \la X,\r\ra )$,


(d) $\X ^c$ is a reversible structure, where $\X ^c =\la X, \r^c\ra$ and $\r ^c =X ^2 \setminus \r$.
\end{fac}
We remark that if $\X$ is a connected binary structure but its complement $\X ^c$ is disconnected  then, by Fact \ref{TB045},
the reversibility of $\X ^c$ implies the reversibility of $\X $.
Thus the results concerning reversibility of disconnected structures can be converted into results about reversibility of such connected structures.
(We note that at least one of the structures $\X$ and $\X ^c$ is connected, see \cite{Ktow})

Reversible structures
have the Cantor-Schr\"{o}der-Bernstein property for condensations. Moreover, for any relational language $L$ we have (see \cite{KuMoVar})
\begin{fac}\label{TA011}
Let $\X$ and $\Y$ be $L$-structures. If $\X$ is a reversible structure and $\Y \sim _c \X$, then
$\Y \cong \X$ (thus $\Y$ is reversible too) and $\Cond (\X ,\Y )=\Iso (\X ,\Y )$.
\end{fac}
\paragraph{Morphisms of disconnected $L_b$-structures}
 If $\X=\la X,\rho  \ra$ is an $L_b$-structure, then
the transitive closure $\sim _\rho$ of the relation $\rho  _{rs} =\Delta _X \cup \rho  \cup \rho  ^{-1}$
(given by $x \sim _\rho y$ iff there are $n\in \N$ and
$z_0 =x , z_1, \dots ,z_n =y$ such that $z_i \;\rho  _{rs} \;z_{i+1}$, for each $i<n$)
is the minimal equivalence relation on $X$ containing $\rho $. The corresponding equivalence classes $[x]$, $x\in X$,
are called the {\it components} of $\X$,
and the structure $\X$ is called {\it connected} iff $|X/\sim _\rho |=1$.

If $\X _i=\la X_i, \rho  _i \ra$, $i\in I$, are connected $L_b$-structures  and $X_i \cap X_j =\emptyset$, for
different $i,j\in I$, then the structure $\bigcup _{i\in I} \X _i =\la \bigcup _{i\in I} X_i , \bigcup _{i\in I} \rho  _i\ra$ is
the {\it disjoint union} of the structures $\X _i$, $i\in I$,  and the structures $\X _i$, $i\in I$, are its components.
A proof of the following fact is direct; see \cite{KuMoSim} for (a) and \cite{Ktow} for (b).
\begin{fac}\label{T4015}
Let $\{\X _i :i\in I\}$ and $\{\Y _j :j\in J\}$ be families of pairwise disjoint and connected
$L_b$-structures and let $\X=\la X,\r\ra$ and $\Y=\la Y,\s\ra$ be their unions. Then

(a) $F\in \Cond (\X ,\Y )$ iff $F=\bigcup _{i\in I} g_i$, where $f:I\stackrel{onto}{\rightarrow}J$,
$g_i \in \Mono (\X _i , \Y _{f(i)})$, for $i\in I$, and
$\{g_i[X_i]: i\in f^{-1}[\{ j \}] \}$ is a partition of $Y_j$, for each $j\in J$;

(b)  $F\in \Emb (\X ,\Y )$ iff $F=\bigcup _{i\in I} g_i$, where $f:I\rightarrow J$, $g_i \in\Emb (\X _i ,\Y _{f(i)})$, for $i\in I$,
and $\la g_i(x), g_{i'}(x')\ra \not\in\s _{rs}$, whenever $i\neq i'$, $x \in X_i$ and $x' \in X_{i'}$.
\end{fac}
\section{Reversibility of disconnected $L_b$-structures}\label{S2}
In this section we prove Theorem \ref{TB065} and give some other equivalents of reversibility in the class of disconnected $L_b$-structures.

\begin{te}\label{TB066}
If $\X _i$, $i\in I$, are pairwise disjoint and connected $L_b$-structures, then
$\bigcup _{i\in I} \X _i$ is reversible iff $\;\bigcup_{i\in J}\mathbb{X}_i$ is reversible for each non-empty set $J\subset I$.

Thus, if $\;\bigcup _{i\in I} \X _i$ is reversible, then all components $\mathbb{X}_i$, $i\in I$, are reversible.
\end{te}
\dok
Let $\X _i =\la X_i,\r_i \ra$, for $i\in I$, and let $\X= \la X,\r\ra=\la \bigcup_{i\in I}X_i,\bigcup_{i\in I}\r_i\ra$.

The implication ``$\Leftarrow$" is trivial.
If there exist a non-empty set $J\subset I$ and $f\in\Cond(\bigcup_{i\in J}\X_i)\setminus \Aut(\bigcup_{i\in J}\X_i)$,
then there are $x,y\in \bigcup_{i\in J}X_i$ such that
$\la x,y \ra \not\in \bigcup_{i\in J}\r _i $ and $\la f(x),f(y) \ra \in \bigcup_{i\in J}\r _i $.
Now, $F:=f\cup \id_{\bigcup_{i\in I\setminus J}X_i}\in\Sym(X)$, it is easy to check
that $F\in \Cond (\X )$,
and the pair $\la x,y\ra$ witnesses that $F\not\in \Aut (\X )$. So $\X$ is not a reversible structure.
\kdok
\begin{te}\label{TA019}
Let $\X _i$, $i\in I$, be pairwise disjoint and connected $L_b$-structures.
Then the structure $\bigcup _{i\in I} \X _i$ is reversible  iff
whenever  $f:I\rightarrow I$ is a surjection,
$g_i \in \Mono(\X _i ,\X _{f(i)})$, for $i\in I$, and
\begin{equation}\label{EQB035}
\forall j\in I \;\; \Big(\Big\{g_i[X_i]: i\in f^{-1}[\{ j\}] \Big\} \mbox{ is a partition of }X_j\Big),
\end{equation}
we have
\begin{equation}\label{EQB034}
f\in \Sym (I) \;\land \;\forall i\in I \;\; g_i \in \Iso (\X _i,\X _{f(i)}) .
\end{equation}
\end{te}
\dok
Let $\X _i =\la X_i,\r_i \ra$, for $i\in I$, and let $\X= \la X,\r\ra=\la \bigcup_{i\in I}X_i,\bigcup_{i\in I}\r_i\ra$.

Suppose that $\X$ is a reversible structure and let the mappings $f$ and $g_i$ be as assumed.
Then by Fact \ref{T4015}(a) we have $F:=\bigcup _{i\in I}g_i\in \Cond (\X )=\Aut (\X )$ and, by  Fact  \ref{T4015}(b),
$g_i \in \Emb (\X _i,\X _{f(i)})$, for all $i\in I$.
Suppose that there are different $i _1,i_2\in I$ such that $f(i_1)=f(i_2)=j$ and let $x_1\in X_{i_1}$ and $x_2\in X_{i_2}$. Since the structure
$\X _j$ is connected there are $y_1, \dots ,y_n \in X_j$ such that $g_{i_1}(x_1)=y_1 (\r _j)_s y_2 \dots (\r _j)_s y_n= g_{i_2}(x_2)$ and, hence, there is $k<n$
such that $y_k$ and $y_{k+1}$ are in different elements of the partition $\{g_i[X_i]: i\in I \land f(i)=j \}$, say
$y_k =g_i (x)\in g_i[X_i] $ and $y_{k+1} =g_{i'} (x')\in g_{i'}[X_{i'}]$, where $i\neq i'$. But then $\la g_i(x), g_{i'}(x')\ra \in\r _{rs}$,
which is, by  Fact \ref{T4015}(b), impossible. Thus $f$ is a bijection and, 
by (\ref{EQB035}), for each $i\in I$ we have $g_i[X_i]=X_{f(i)}$ and,
hence, $g_i \in \Iso (\X _i,\X _{f(i)})$.

Conversely, for $F\in \Cond (\X )$ we prove that $F\in \Aut (\X )$. By  Fact  \ref{T4015}(a) and the assumption we have
$F=\bigcup _{i\in I} g_i$, where $f\in \Sym (I)$ and $g_i \in \Iso (\X _i,\X _{f(i)})$, for all $i\in I$. By  Fact  \ref{T4015}(b)
we have $F\in \Emb (\X )$ and, since $F$ is a surjection, $F\in \Aut (\X )$.
\kdok
Now we prove Theorem \ref{TB065}. In fact we will prove its contrapositive.
\begin{te}\label{TB014}
($\ast$) The union $\X:=\bigcup _{i\in I}\X_i$ is not reversible iff $\la \X _i :i\in I\ra$ is not a reversible sequence of structures or
there is a non-trivial $\Z$-sequence in $I$.
\end{te}
\dok
$(\Rightarrow)$
If $\X$ is not a reversible structure, then
by Theorem \ref{TA019}
there are $f\in \Sur (I)$ and $g_i\in\Mono(\mathbb{X}_i,\mathbb{X}_{f(i)})$, for $i\in I$, satisfying (\ref{EQB035}) and $\neg \,$(\ref{EQB034}).
By (\ref{EQB035}) and Fact \ref{T4015}(a) we have
$\bigcup_{i\in f^{-1}[\{j\}]}g_i\in\Cond (\bigcup_{i\in f^{-1}[\{j\}]}\mathbb{X}_i,\mathbb{X}_j)$, for all $j\in I$.
So, if $f\not\in \Sym (I)$, then we have $\neg \,$(\ref{EQB040}). 

If $f\in \Sym (I)$, then, by $\neg\,$(\ref{EQB034}),  $g_{i_0}\not\in\Iso(\mathbb{X}_{i_0},\mathbb{X}_{f(i_0)})$, for some $i_0\in I$.
Since $g_{i_0}\in\Cond(\mathbb{X}_{i_0},\mathbb{X}_{f(i_0)})$
and, since by the reversibility of $\mathbb{X}_{i_0}$ 
and Fact \ref{TA011}, $\mathbb{X}_{i_0}\cong\mathbb{X}_{f(i_0)}$ would imply
$\Cond(\mathbb{X}_{i_0},\mathbb{X}_{f(i_0)})=\Iso(\mathbb{X}_{i_0},\mathbb{X}_{f(i_0)})$, we have $\mathbb{X}_{i_0}\not\cong\mathbb{X}_{f(i_0)}$.

Let $i_k:=f^{k}(i_0)$, for $k\in\mathbb{Z}$. Then $\mathbb{X}_{i_0}\not\cong\mathbb{X}_{i_1}$ and
for each $k\in\mathbb{Z}$ we have $f^{-1}[\{ i_{k+1}\}]=\{ i_{k}\}$ and, hence, $g_{i_{k}}\in\Cond(\mathbb{X}_{i_{k}},\mathbb{X}_{i_{k+1}})$.

Suppose  $i_k=i_l$, for some $k<l\in \Z$.
Then $f^{l}(i_0)=f^{k}(i_0)$, and, hence, $i_{l-k}=f^{l-k}(i_0)=i_0$, where $l-k\geq 1$. 
So,
$\mathbb{X}_{i_0}\preccurlyeq _c\mathbb{X}_{i_1}\preccurlyeq _c\cdots\preccurlyeq _c\mathbb{X}_{i_{l-k}}=\mathbb{X}_{i_0}$,
which implies $\mathbb{X}_{i_0}\sim_c\mathbb{X}_{i_1}$.
Since $\mathbb{X}_{i_0}$ is a reversible structure,  by Fact \ref{TA011} we would have $\mathbb{X}_{i_0}\cong\mathbb{X}_{i_1}$, which is false.
So $\la i_k \!:\!k\in \Z\ra$ is an injection and, since $\mathbb{X}_{i_0}\!\not\cong\!\mathbb{X}_{i_1}$, it is a non-trivial $\Z$-sequence in $I$.

($\Leftarrow$) If there are $f\in \Sur(I)\setminus \Sym (I)$ and $G_j\in\Cond(\bigcup_{i\in f^{-1}[\{j\}]}\mathbb{X}_i,\mathbb{X}_j)$,
for $j\in I$, then, clearly,
$g_i:=G_{f(i)}\!\upharpoonright X_i\in\Mono(\mathbb{X}_i,\mathbb{X}_{f(i)})$, for all $i\in I$,
and (\ref{EQB035}) is true.
Since $f\not\in \Sym (I)$, by Theorem \ref{TA019} the structure $\X$ is not reversible.

Suppose that $\langle i_k:k\in\mathbb{Z}\rangle$ is a $\Z$-sequence in $I$ and that $\X_{i_r}\not\cong\X_{i_{r+1}}$, for some $r\in \Z$.
Then the function $f:I\rightarrow I$, defined by $f(i)=i$, for $i\in I\setminus\{i_k:k\in\mathbb{Z}\}$,
and $f(i_k)=i_{k+1}$, for $k\in\mathbb{Z}$, is a bijection.
Let $g_i:=\id_{X_i}$, for $i\in I\setminus\{i_k:k\in\mathbb{Z}\}$,
and let us take $g_{i_k}\in\Cond(\mathbb{X}_{i_k},\mathbb{X}_{i_{k+1}})$, for $k\in\mathbb{Z}$.
Then $g_i\in\Mono(\mathbb{X}_i,\mathbb{X}_{f(i)})$, for all $i\in I$, and (\ref{EQB035}) holds.
But $g_{i_r}\not\in\Iso(\mathbb{X}_{i_r},\mathbb{X}_{i_{r+1}})$ and, by Theorem \ref{TA019}, the structure $\X$ is not reversible.
\kdok
\begin{cor}\label{TB051}
An $L_b$-structure with finitely many components is reversible
iff all its components are reversible.
\end{cor}
\dok
Let $\X =\bigcup _{i\in I}\X _i$, where $|I|<\o$ and $\X _i$, $i\in I$, are pairwise disjoint and connected $L_b$-structures.
The implication ``$\Rightarrow$" follows from Theorem \ref{TB066}. If the structures $\X _i$, $i\in I$, are reversible,
then ($\ast$) holds, so, since $\Sur (I)=\Sym (I)$ and there are no $\Z$-sequences in $I$, by Theorem \ref{TB014} the structure $\X$ is reversible.
\kdok
\begin{cor}\label{TB064}
An $L_b$-structure $\bigcup _{i\in I}\X _i$ with finite components is reversible
iff $\la \X _i :i\in I\ra$ is a reversible sequence of structures
and there are no infinite classes $[ \X _i]_{\cong },[ \X _j]_{\cong }\in \CX /\!\cong$ such that $\X _i \prec _c \X _j $.
\end{cor}
\dok
Since all finite structures are reversible condition ($\ast$) is fulfilled.
According to Theorem \ref{TB065} we show that the negation of the second condition holds iff there is a non-trivial $\Z$-sequence in $I$.
The implication ``$\Rightarrow$" is trivial.

If $\la i_k :k\in \Z \ra$ is a non-trivial $\Z$-sequence in $I$, then $\X _{i_k} \preccurlyeq _c \X _{i_{k+1}}$, for all $k\in \Z$,
and there is $k_0\in \Z$ such that $\X _{i_{k_0}} \prec _c \X _{i_{k_0+1}}$.
Since $\X _{i_{k_0+1}}\preccurlyeq _c \X _{i_{k_0+2}}\preccurlyeq _c \dots$ and the structures $\X _i$, $i\in I$, are finite,
there is $s\geq k_0+1$ such that $\X _{i_k}\cong\X _{i_s}$, for all $k\geq s$, and, similarly,
there is $r\leq k_0$ such that $\X _{i_k}\cong\X _{i_r}$, for all $k\leq r$.
Now we have $\X _{i_{r}} \prec _c \X _{i_{s}}$ and the classes $[\X _{i_{r}}]_{\cong }$ and $[\X _{i_{s}}]_{\cong }$ are infinite.
\kdok
By Theorem 3.4 of \cite{KuMoEqu}, if $\X _i$, $i\in I$, are pairwise disjoint tournaments (resp.\ in particular, linear orders),
and  $\la |X_i|:i\in I\ra$ is a reversible sequence of cardinals, then the digraph (resp.\ poset) $\bigcup _{i\in I}\X_i$ is reversible.
But that condition is not necessary for the reversibility of such unions and now we give a characterization.
\begin{cor}\label{TB052}
A disjoint union $\bigcup _{i\in I}\X_i$  of linear orders (or, more generally, tournaments) is not reversible iff there is a
non-injective surjection $f:I\rightarrow I$ such that each component $\X _j$ can be partitioned into copies of $\X _i$,
where $i\in f^{-1}[\{ j\}]$.\footnote{that is, there is a
partition $\{ A_i: i\in f^{-1}[\{ j\}]\}$ of $X_j$
such that $\A_i \cong \X _i$, for all $i\in f^{-1}[\{ j\}]$.}
\end{cor}
\dok
($\ast$) is true, because all tournaments are reversible and connected. Since  for any two tournaments
$\X$ and $\Y$ we have $\Cond(\X ,\Y )= \Iso(\X ,\Y)$, all $\Z$-sequences in $I$ are trivial and, by Theorem \ref{TB014},
$\bigcup _{i\in I}\X_i$ is not reversible iff there is $f\in \Sur (I)\setminus \Sym (I)$ such that for each $j\in I$
there is $G_j \in \Cond (\bigcup_{i\in f^{-1}[\{ j\}]}\X _i, \X _j)$.

Then for $i\in f^{-1}[\{ j\}]$ and $A_i:=G_j[X_i]$ we have
$G_j\!\upharpoonright \!X_i \in \Cond(\X_i , \A_i)=\Iso(\X_i , \A_i)$ and $\{ A_i: i\in f^{-1}[\{ j\}]\}$
is a partition of $X_j$.

Conversely, if $f\in \Sur (I)\setminus \Sym (I)$ and
$\{ A_i: i\in f^{-1}[\{ j\}]\}$ is a partition of $X_j$, for each $j\in I$, and $g_i \in \Iso(\X_i , \A_i)$, for $i\in f^{-1}[\{ j\}]$,
then, clearly, $G_j:=\bigcup _{i\in f^{-1}[\{ j\}]\}}g_i\in \Cond (\bigcup_{i\in f^{-1}[\{ j\}]}\X _i, \X _j)$
and we are done.
\kdok
\section{Triviality of $\o^*$-sequences of monomorphisms}\label{S3}
In this and the following section we consider some conditions which imply reversibility of disconnected $L_\b$-structures.
If $\X$ and $\Y$ are $L_b$-structures and there is a monomorphism $f:\X \rightarrow \Y$, we will write $\X \preccurlyeq_m \Y$.

Under ($\ast$), a mapping $i:\o \rightarrow I$, usually denoted by $\la i_k :k\in \o\ra$, will be called an
{\it $\o ^*$-sequence in $I$} iff it is an injection and
\begin{equation}\label{EQB042}\textstyle
\forall k,l\in \o \;\;( k<l \Rightarrow \X_{i_l}\preccurlyeq_m \X_{i_k}).
\end{equation}
If, in addition, $\Mono(\X _{i_1},\X _{i_0})=\Iso(\X _{i_1},\X _{i_0})$,
the $\o ^*$-sequence $\la i_k :k\in \o\ra$ will be called {\it trivial}.
By the transitivity of the relation $\preccurlyeq_m$, condition
(\ref{EQB042}) is equivalent to the existence of a sequence of monomorphisms
$g_k:\X_{i_{k+1}}\rightarrow\X_{i_k}$, $k\in \o$.\footnote{We note that $\la i_k :k\in \o\ra$ is an $\o ^*$-sequence in $I$ iff
$k\mapsto \X _{i_k}$ is a monomorphism from the linear order $\o^*=\la \o ,>\ra$ to the preorder $\la \CX ,\preccurlyeq_m\ra$
and $\la i_k :k\in \o\ra$ is non-trivial iff we can choose $g_0:\X_{i_{1}}\rightarrow\X_{i_0}$ which is not an isomorphism.
This holds if, in particular, $\X_{i_{1}}\not\cong \X_{i_0}$.}
\begin{fac}\label{TB048}
Let $f:I\rightarrow I$ be a surjection, let $j\in I$, where $|f^{-1}[\{j\}|>1$, and let $O(j):=\{f^n(j):n\in\omega\} $. Then
$|f^{-1}[\{j\}]\cap O(j)|\leq 1$.
\end{fac}
\dok
If there exists $i\in f^{-1}[\{j\}]\cap O(j)$, then $i=f^k(j)$, for some $k\in \o$, and  $j=f(i)= f^{k+1}(j)$.
So, there is $l =  \min \{ m\in\omega : f^m(j)= j\}$.

If  $l=0$, that is $f(j)=j$, then $O(j)=\{j\}$ and the statement is true.

If $l>0$, then $f(f^{l-1}(j))=f^l(j)=j$ and, hence, $f^{l-1}(j)\in f^{-1}[\{j\}]\cap O(j)$.
Clearly, $f^{ql}(j)=j$, for all $q\in \o$.
So, if $n=ql +r$, where $q\in \o$ and $r<l$, then $f^{n}(j)=f^r (f^{ql}(j))=f^r (j)$
and, hence, $O(j)=\{ f^{n}(j): n\leq l-1\}$. Assuming  that  $f^{n}(j)\in f^{-1}[\{j\}]\cap O(j)$, for some $n<l-1$,
we would have $f^{n+1}(j)=j$ and $n+1<l$, which contradicts the minimality of $l$.
Thus  $f^{-1}[\{j\}]\cap O(j)=\{f^{l-1}(j)  \}$ and we are done.
\kdok
\begin{te}\label{TB047}
($\ast$) If each $\o ^*$-sequence in $I$ is trivial,\footnote{Then each monomorphism $i:\o^* \rightarrow\la \CX ,\preccurlyeq_m\ra$
maps $\o$ into the $\cong$-class of $\X _{i_0}$ (if $\X _{i_0}\not\cong \X _{i_k}$, then $\la i_0, i_k, i_{k+1}, \dots\ra$ is a non-trivial $\o ^*$-sequence in $I$)
but the converse is not true (take $\bigcup _\o \o$).} the structure $\X$ is reversible.
\end{te}
\dok
If $\X$ is not reversible, then by Theorem \ref{TB014} we have the following two cases.

{\it Case 1.} There is a $\Z$-sequence $\la i_k:k\in\Z\ra$ in $I$ such that $\X_{i_{-1}}\not\cong\X_{i_0}$. Let $j_k:=i_{-k}$, for $k\in\mathbb{Z}$. Then for each $k\in\omega$ we have
$\Mono(\mathbb{X}_{j_{k+1}},\mathbb{X}_{j_k})\supset\Cond(\mathbb{X}_{i_{-(k+1)}},\mathbb{X}_{i_{-k}})\neq\emptyset$ and, hence,
$\langle j_k:k\in\omega\rangle$ is an  $\o ^*$-sequence in $I$
and $\Mono(\mathbb{X}_{j_1},\mathbb{X}_{j_0})\neq\Iso(\mathbb{X}_{j_1},\mathbb{X}_{j_0})=\emptyset$.

{\it Case 2.} There is $f\in \Sur(I)\setminus \Sym (I)$ such that for each $j\in I$ there exists
$G_j\in\Cond(\bigcup_{i\in f^{-1}[\{j\}]}\mathbb{X}_i,\mathbb{X}_j)$.
Let $j_*\in I$, where $|f^{-1}[\{j_*\}|>1$. By Fact \ref{TB048} there is
\begin{equation}\label{EQB038}
i_*\in f^{-1}[\{j_*\}]\setminus\{f^n(j_*):n\in\omega\}.
\end{equation}
Since $f:I\rightarrow I$ is a surjection, there is a sequence $\langle i_k:k\in\omega\rangle\in I^\omega$, such that
$i_0=j_*$, $i_1=i_*$, and $f(i_{k+1})=i_k$, for all $k\in\omega$. Suppose that there is $k\in \N$ such that
$i_k \in \{f^n(j_*):n\in\omega\}$ and let $k$ be the minimal such element of $\N$.
By (\ref{EQB038}) we have $k>1$. So $i_k = f^n(j_*)$, for some $n\in \o$, and, hence, $i_{k-1}=f(i_k)=f^{n+1}(j_*)$,
which contradicts the minimality of $k$. Thus
\begin{equation}\label{EQB039}
\{i_k:k\in\mathbb{N}\}\cap\{f^n(j_*):n\in\omega\}=\emptyset .
\end{equation}
Suppose that $i:\omega\rightarrow I$ is not an injection and let $r$ be the minimal element of $\omega$ such that $i_r=i_s$, for some $s>r$.
$r=0$ would imply that $i_s=i_0 =j_*\in \{f^n(j_*):n\in\omega\}$, which is impossible by (\ref{EQB039}).
Now $i_{r-1}=f(i_r)=f(i_s)=i_{s-1}$, which is impossible by the minimality of $r$.
Thus $i:\omega\rightarrow I$ is an injection.

For $k\in\omega$ we have $G_{i_k}\in\Cond(\bigcup_{i\in f^{-1}[\{i_k\}]}\mathbb{X}_i,\mathbb{X}_{i_k})$ and $i_{k+1}\in f^{-1}[\{i_k\}]$
so $G_{i_k}\!\upharpoonright X_{i_{k+1}}\in\Mono(\mathbb{X}_{i_{k+1}},\mathbb{X}_{i_k})$.
Since $|f^{-1}[\{i_0\}]|>1$, is follows that
$G_{i_0}\!\upharpoonright_{\,X_{i_1}}\in\Mono(\mathbb{X}_{i_1},\mathbb{X}_{i_0})\setminus\Iso(\mathbb{X}_{i_1},\mathbb{X}_{i_0})$. Therefore, $\langle i_k:k\in\omega\rangle$ is a non-trivial $\omega^*$-sequence in $I$.
\kdok
\begin{ex}\label{EXB009}\rm
The converse of Theorem \ref{TB047} is not true.
The equivalence relation $\bigcup _{i\in \o}\X _i$, where $|X_i|=2$, for even $i$'s, and $|X_i|=5$, for odd $i$'s,
is reversible \cite{KuMoEqu} (see, also, Introduction) and $\la 1,0,2,4,6,\dots\ra$ is a non-trivial $\o ^*$-sequence in $\o$.
\end{ex}
\begin{rem}\label{RB000}\rm
If each $\o ^*$-sequence in $I$ is trivial, then we have two possibilities:

1. There is an $\o ^*$-sequence $\la i_k:k\in \o\ra$. Then there is
$g_0\in\Mono(\X _{i_1},\X _{i_0})=\Iso(\X _{i_1},\X _{i_0})$, which implies that $\X _{i_1}\cong\X _{i_0}$ and
$\Mono(\X _{i_0})=\Aut(\X _{i_0})$. Generally speaking, if $\X$ is a structure satisfying $\Mono(\X )=\Aut(\X )$, then
$\Cond(\X )=\Emb(\X )=\Aut(\X )$; so it is reversible and copy-minimal (see \cite{Kmin} for examples). Clearly, if $\X$ is a finite structure,
then $\Mono(\X )=\Aut(\X )$ and the linear graph $\BG _\Z$ is a reversible, connected infinite structure satisfying $\Mono(\X )=\Aut(\X )$.
We note that in \cite{Dus} Dushnik and  Miller constructed {\it embedding-rigid} dense suborders $\BL$ of the real line
(i.e.\ $\Emb (\BL)=\{ \id _L\}$; see also \cite{Rosen}, p.\ 147) of size $\c$ and
similar examples can be made using the ZFC result of  Vop\v enka,  Pultr and Hedrl\'in \cite{Vop}
saying that on every set there is an endo-rigid binary relation.

2. $\o ^*$-sequences do not exist at all. This situation is considered in the sequel.
\end{rem}
\section{Non-existence of $\o ^*$-sequences. Monotone functions}\label{S4}
We recall that a pair $\CW =\la \CA ,\CR \ra$ is called a {\it well founded relation} (we will write $\CW \in \wfr$) iff
$\CA$ is a class, $\CR$ a binary class relation on $\CA$ and each non-empty {\it set} $X\subset \CA$
has an $\CR$-minimal element, that is,
\begin{equation}\label{EQB043}
\forall X \;\; \Big(\emptyset \neq X \subset \CA \Rightarrow  \exists y\in X \;\neg \exists z\in X \;\; z\CR y \Big).\footnote{For
convenience, we use that general notions from set theory. The classes  $\CA$ and $\CR$ are collections of sets satisfying some formulas,
say $A(v)$ and $R(u,v)$ of the language of set theory, $\{ \in \}$. For example, $\CA$ is the class of all ordinals and $\CR$ the usual strict linear order on that class.}
\end{equation}
Note that then the relation $\CR$ on $\CA$ is irreflexive and asymmetric (a ``class-digraph") and its reflexivization $\leq _{\CR }$ is defined  by
$$
a\leq _{\CR } b \;\Leftrightarrow\; a=b \;\lor \;a\CR b .
$$
If, in addition,  $\CC \subset \Mod _L$ is a class of $L$-structures, we will say that a (class) function $\theta : \CC \rightarrow \CA$
is {\it monotone with respect to monomorphisms} iff
\begin{equation}\label{EQB044}
\forall \X ,\Y \in \CC \;\; \Big(\X  \preccurlyeq _m \Y \Rightarrow  \theta (\X ) \leq _{\CR } \theta (\Y ) \Big)
\end{equation}
and the class of such functions (which are, in fact, homomorphisms from the preorder
$\la \CC , \preccurlyeq _m\ra$ to the reflexivization of $\CW$) will be denoted by $\M (\CC ,\CW)$.

If ($\ast$) holds, let $\M (\CX )=\bigcup _{\CW \in \wfr}\M (\CX , \CW)$.
For $\theta \in \M (\CX ,\CW)$ and $a \in \CA$, let
\begin{equation}\label{EQB049}
I^\theta _a := \{ i\in I : \theta (\X _i )=a \}.
\end{equation}
\begin{te}\label{TB050}
($\ast$) If there is $\theta \in \M (\CX )$ such that for each $a\in \theta [\CX ]$ there are no $\o ^*$-sequences in $I^\theta _a$,
then $\X$ is a reversible structure.
This holds if, in particular, $\la \theta (\X _i) :i\in I\ra$ is a finite-to-one sequence
(i.e., if the sets $I^\theta _a$, $a\in \theta [\CX ]$, are finite).
\end{te}
\dok
First we prove that the following conditions are equivalent:

(i) There is an $\o ^*$-sequence in $I$,

(ii) $\forall \theta \in \M (\CX )\;\;
     \exists a\in \theta [\CX ] \;\; \exists \;\o ^*\mbox{-sequence in }I^\theta _a$,

(iii) $\exists \theta \in \M (\CX )\;\;
     \exists a\in \theta [\CX ] \;\; \exists \;\o ^*\mbox{-sequence in }I^\theta _a$.

\noindent
(i) $\Rightarrow$ (ii). Let $\la i_k :k\in \o\ra$ be an $\o ^*$-sequence in $I$, $\CW \in \wfr$ and $\theta \in \M (\CX ,\CW)$.
Since $X=\{ \theta (\X _{i_k}): k\in \o \}$ is a non-empty subset of $\CA$, by (\ref{EQB043}) there is $k_0\in \o$
such that
\begin{equation}\label{EQB045}
\forall k\in \o \;\; \neg \,\theta (\X _{i_k})\CR \theta (\X _{i_{k_0}}) .
\end{equation}
Now  for $k\geq k_0$ we have $\X _{i_k}\preccurlyeq _m \X _{i_{k_0}}$, which
by (\ref{EQB044}) implies $\theta (\X _{i_k})\leq _{\CR } \theta (\X _{i_{k_0}}) =: a$ and, by (\ref{EQB045}), $\theta (\X _{i_k})= \theta (\X _{i_{k_0}})$.
So $\theta (\X _{i_k})= a$, for all $k\geq k_0$, that is $\{ i_k : k\geq k_0 \}\subset I^\theta _a$ and $\la i_{k_0 +k}:k\in \o\ra$
is an $\o ^*$-sequence in $I^\theta _a$.

(ii) $\Rightarrow$ (iii). If $\Card$ denotes the class of all cardinals, then, clearly, $\CW =\la \Card, <\ra\in \wfr$
and for the function $\theta :\CX \rightarrow \Card$ defined by $\theta (\X _i)=|X_i|$ we have $\theta \in \M (\CX ,\CW)$.
By (ii) there are a cardinal $\k$ and an $\o ^*$-sequence in $I^\theta _\k$.

(iii) $\Rightarrow$ (i). This is trivial, since an $\o ^*$-sequence in $I^\theta _a$ is an $\o ^*$-sequence in $I$.

\noindent
Now, by the assumption, $\neg \,$(ii) holds and, hence, there are no $\o ^*$-sequences in $I$ and we apply Theorem \ref{TB047}.
\kdok
\begin{cor}\label{TB056}
If {\rm ($\ast$)} holds and the sequence
$\la |X_i|:i\in I\ra$ is finite-to-one,\footnote{that is,
there is no infinite $J\subset I$ such that $|\X _i| = |\X _j|$, for all $i,j\in J$.} then
$\bigcup _{i\in I}\X_i$ is a reversible structure.
\end{cor}
For the structures with finite components the condition implying the reversibility of $\bigcup _{i\in I}\X_i$
given in Theorem \ref{TB050} is, in fact, equivalent to a simpler condition.
\begin{prop}\label{TB063}
If ($\ast$) holds and the structures $\X _i$, $i\in I$, are finite, then the following conditions are equivalent:

(a) The sequence $\la |X _i|:i\in I\ra$ is finite-to-one,

(b) $\exists \theta \in \M (\CX )\;\; \forall a\in \theta [\CX ]  \;\;\neg \exists \;\o ^*$-sequence in $I^\theta _a$.
\end{prop}
\dok
(a) $\Rightarrow$ (b) is trivial: take $\theta : \CX \rightarrow \Card$, where $\theta (\X _i)=|X_i|$.

$\neg\,$(a) $\Rightarrow$ $\neg\,$(b). Let $J\subset I$ and $n\in \N$, where $|J|\geq \o$ and $|X_i|=n$, for all $i\in J$.
Then, since the structures $\X _i$, $i\in I$, are finite, there is $K\subset J$, where $|K|\geq \o$ and $\X _i\cong \X _j$, for all $i,j\in K$;
let us fix an $i_0\in K$. Let $\theta \in \M (\CX, \CW )$ and $a\in \theta [\CX ]$, where $\theta (\X _{i_0})=a$.
Now, if $i\in K$, then $\X_i \cong \X_{i_0}$ and, hence $\X_i \preccurlyeq _m \X_{i_0} \preccurlyeq _m \X_i $, which implies that
$\theta (\X_i) \leq _{\CR} \theta (\X_{i_0}) \leq _{\CR} \theta (\X_i )$, and, thus, $\theta (\X_i )=a$, that is
$i\in I^\theta _a$. So $K\subset I^{\theta } _a$ and taking an injection $i: \o \rightarrow K$ we obtain
an $\o ^*$-sequence $\la i_k : k\in \o\ra$ in $I^\theta _a$
(because $\X _{i_{k+1}} \cong \X _{i_{k}}$ gives $\X _{i_{k+1}} \preccurlyeq _m \X _{i_{k}}$).
\kdok
\paragraph{Finite diagonal products of monotone functions}
The class $\wfr$ is not closed under direct products (in the product of two-element chains, $2^\o$,
the set $X=\{ x_n:n\in \o\}$, where $x_n =\la 0, \dots,0,1,1,\dots\ra$ has $n$-many zeros, does not have a minimal element).
But $\wfr$ is closed under finite products.
\begin{te}\label{TB057}
Let $n\in \N$ and let $\CW _k =\la \CA _k,\CR _k\ra\in \wfr$, for $k<n$. Then

(a) $\la \prod _{k<n} \CA _k,\CR\ra\in \wfr$, where for $a=\la a_k \ra$, $b=\la b_k \ra\in \prod _{k<n} \CA _k$ we have
\begin{equation}\label{EQB046}
a \CR b \Leftrightarrow \forall k<n \;\;( a_k = b_k \lor a_k \CR _k b_k) \land \exists k<n \;\; a_k \CR _k b_k;
\end{equation}

(b) If $\CC\subset \Mod _L$ and $\theta _k \in \M (\CC , \CW _k )  $, for $k<n$,
then $\theta \in \M (\CC , \prod _{k<n} \CA _k )$,
where $\theta$ is the diagonal mapping $\theta =\Delta _{k<n}\theta _k : \CC \rightarrow \prod _{k<n} \CA _k$, defined by
\begin{equation}\label{EQB047}
\theta (\X )=\la \theta _k (\X ): k<n\ra ;
\end{equation}

(c) If ($\ast$) holds and in (b) we put $\CC=\CX$, then for $a=\la a_k :k<n \ra\in \theta [\CX ]$ we have
\begin{equation}\label{EQB048}\textstyle
I^\theta _a= \{ i\in I : \forall k<n \;\;\theta _k (\X )=a_k\}=\bigcap _{k<n}I^{\theta _k}_{a_k}.
\end{equation}
Thus the partition $\{ I^\theta _a: a\in \theta [\CX ]\}$ of $I$ refines all partitions $\{ I^{\theta _k} _{a_k}: a_k\in \theta _k[\CX ]\}$, $k<n$,
and, hence, the sequence $\la \theta (\X _i): i\in I \ra$ has more chance to be finite-to-one (see Theorem \ref{TB050}).
\end{te}
\dok
(a) Suppose that a non-empty set $X \subset \prod _{k<n} \CA _k$ has no $\CR$-minimal elements.
Then there are $a^r \in X$, $r\in \o$, such that for each  $r\in \o$ we have $a^{r+1}\CR a^r$ and, by  (\ref{EQB046}),
there is $k_r <n$ satisfying $a^{r+1}_{k_r}\CR_{k_r} a^r_{k_r}$. Thus there are $k^*<n$ and an increasing sequence $\la r_s :s\in \o\ra$ in $\o$
such that
\begin{equation}\label{EQB050}\textstyle
\forall s\in \o \;\; (k_{r_s}=k^* \;\land \;a^{r_s+1}_{k^*}\CR _{k^*}a^{r_s}_{k^*})
\end{equation}
\begin{equation}\label{EQB051}\textstyle
\forall r\in \o \setminus  \{ r_s :s\in \o\} \;\; a^{r+1}_{k^*}= a^{r}_{k^*}.
\end{equation}
If $r\in \o$ and $s^*:=\min \{ s\in \o : r_s\geq r\}$, then $r_{s^*}\geq r$ and by (\ref{EQB050}) and (\ref{EQB051}) we have
$a^{r_{s^*}+1}_{k^*}\CR _{k^*}a^{r_{s^*}}_{k^*}
=a^{r_{s^*}-1}_{k^*}
=a^{r_{s^*}-2}_{k^*}
=\dots
=a^{r}_{k^*}$.
This implies that the subset $\{ a^{r}_{k^*}: r\in \o\}$ of $\CA _{k^*}$ has no $\CR_{k^*}$-minimal elements, which
contradicts the assumption that $\CW _{k^*} \in \wfr$.

(b) If $\X ,\Y \in \CC$ and $\X \preccurlyeq _m \Y$, then, by the assumption,
for each $k<n$ we have $\theta _k (\X )= \theta _k (\Y )$ or $\theta _k (\X ) \CR _k \theta _k (\Y )$.
So, if $\theta _k (\X )= \theta _k (\Y )$, for all $k<n$, then by (\ref{EQB047}) we have $\theta  (\X )= \theta  (\Y )$.
Otherwise, there is $k<n$ such that $\theta _k (\X ) \CR _k \theta _k (\Y )$ and, by (\ref{EQB046}), $\theta _k (\X )\CR \theta _k (\Y )$.
Statement (c) follows from (\ref{EQB047}) and (\ref{EQB049}).
\kdok
\paragraph{Some examples of diagonal products}
Let $\LO$, $\PO$, $\Ord$ and $\Ord ^*$ denote the classes of linear orders, partial orders, ordinals and reversed ordinals respectively.
\begin{fac}\label{TB061}
Let $\X ,\Y \in \PO$.

(a) If $\BL \in \LO$, then $\Mono (\BL ,\X )=\Emb (\BL ,\X )$;

(b) $\{ \BL \in \LO : \BL \preccurlyeq _m \X \}= \{ \BL \in \LO : \BL \hookrightarrow \X \}$;

(c) If $\X \preccurlyeq _m \Y$, then $\{ \BL \in \LO : \BL \hookrightarrow \X\}\subset \{ \BL \in \LO : \BL \hookrightarrow \Y\}$.
\end{fac}
\dok
(a) Let $f: \BL \rightarrow \X$ be a monomorphism. If $x,y\in L$ and $f(x)<_\X f(y)$, then $x\neq y$ and,
since $f$ is a homomorphism, $y<x$ would imply $f(y)<_\X f(x)$, which is not true.
So, since $\BL$ is a linear order, $x<y$, and, thus, $f$ is a strong homomorphism.
Since $f$ is one-to-one it is an embedding. (b) follows from (a).

(c) Let $f: \X \rightarrow \Y$ be a monomorphism and let $g:\BL \hookrightarrow \X$. Then $f\upharpoonright g[L]$
is a monomorphism from the linear order $\la g[L], <_\X \upharpoonright g[L] \ra$ into the poset $\Y$ and, by (a)
it is an embedding. Thus $(f\upharpoonright g[L]) \circ g:\BL \hookrightarrow \Y$.
\kdok
Let the (class) functions $\theta _0,\theta _1: \PO \rightarrow \Ord$ be defined by:
$$\theta _0 (\X )= \sup \{ \a \in \Ord : \a \hookrightarrow \X\} \mbox{ and }
\theta _1 (\X )= \sup \{ \a \in \Ord : \a ^* \hookrightarrow \X\}.$$
\begin{prop}\label{TB062}
If {\rm ($\ast$)} holds, $\X_i\in \PO$, for $i\in I$, and
$\la \la \theta _0 (\X _i),\theta _1 (\X _i)\ra :i\in I\ra$ is a finite-to-one sequence, then
$\bigcup _{i\in I}\X_i$ is a reversible poset.
\end{prop}
\dok
If $\X ,\Y \in \PO$ and $\X \preccurlyeq _m \Y$, then by Fact \ref{TB061}(c)
$\{  \a \in \Ord : \a \hookrightarrow \X\}\subset \{ \a \in \Ord : \a \hookrightarrow \Y\}$
and, hence, $\theta _0 (\X )\leq \theta _0(\Y )$. So $\theta _0\in \M (\PO ,\Ord )$ and, similarly,
$\theta _1\in \M (\PO ,\Ord )$. By Theorem \ref{TB057}(b) we have $\theta \in \M(\PO , \Ord \times \Ord )$,
where $\theta : \PO \rightarrow \Ord \times \Ord$ is defined by  $\theta (\X )=\la \theta _0 (\X ),\theta _1 (\X )\ra$.
Now the statement follows from Theorem  \ref{TB050}.
\hfill $\Box$
\begin{ex}\label{EXB007}\rm
Let $I$ be the set of pairs of countably infinite ordinals, that is $I=(\o _1 \setminus \o)^2$,
and let $\X_{\la \a ,\b\ra}$, for $\la \a ,\b\ra \in I$, be disjoint partial orders
such that, using the notation from Proposition \ref{TB062},
$\theta _0 (\X _{\la \a ,\b\ra})=\a$ and $\theta _1 (\X _{\la \a ,\b\ra})=\b$. Then
$\bigcup _{\la \a ,\b\ra \in I}\X _{\la \a ,\b\ra}$ is a reversible poset.
If, in particular, $X_{\la \a ,\b\ra}\cong \b ^* +\a$, this follows from Corollary \ref{TB055} as well.
We note that here the sequences $\la \theta _0 (\X _i):i\in I\ra$ and $\la \theta _1 (\X _i):i\in I\ra$
are not finite-to-one, but $\la \theta (\X _i):i\in I\ra$ is one-to-one.
\end{ex}
\section{Applications: disjoint unions of chains}\label{S5}
\paragraph{$\s$-scattered chains}
We recall that a linear order (chain) $\BL$ is called {\it scattered}, we will write $\BL\in \Scatt$, if it does not contain a dense suborder
(equivalently, iff $\Q \not\hookrightarrow \BL$);
$\BL$ is said to be {\it $\s$-scattered}, we will write $\BL\in \SScatt$, iff $\BL$ is at most countable union of scattered linear orders.
\begin{prop}\label{TB053}
If {\rm ($\ast$)} holds and $\X_i$,  $i\in I$, are $\s$-scattered linear orders, then\\[-2mm]

$\la [\X _i]_\rightleftarrows :i\in I\ra$ is a finite-to-one sequence\footnote{that is,
there is no infinite $J\subset I$ such that $\X _i \rightleftarrows \X _j$, for all $i,j\in J$.} $\Rightarrow$
$\bigcup _{i\in I}\X_i$ is a reversible poset.
\end{prop}
\dok
Clearly, $\la \SScatt ,\hookrightarrow \ra$ is a preorder,
the bi-embedability relation $\rightleftarrows$ defined on $\SScatt$ by:
$\BL \rightleftarrows \BL ' \Leftrightarrow \BL \hookrightarrow \BL ' \land \BL ' \hookrightarrow \BL $, is an equivalence relation,
and, denoting the equivalence class of $\BL$ by $[\BL ]_\rightleftarrows$,
we obtain the corresponding antisymmetric quotient, $\la \SScatt /\rightleftarrows ,\trianglelefteq \ra$,
where the partial order $\trianglelefteq$ is defined by
$[\BL ]_\rightleftarrows \trianglelefteq [\BL ']_\rightleftarrows \Leftrightarrow \BL \hookrightarrow \BL ' $.
Writing $\BL \prec \BL '$ iff $\BL \hookrightarrow \BL ' \land \BL ' \not\hookrightarrow \BL$,
the corresponding strict (irreflexive) partial order is the structure $\la \SScatt /\rightleftarrows ,\vartriangleleft \ra$, where
$$
[\BL ]_\rightleftarrows \vartriangleleft [\BL ']_\rightleftarrows
\;\Leftrightarrow \; \BL \prec \BL '.
$$
From the classical Laver's result (that $\la \SScatt , \hookrightarrow \ra$ is a better-quasi-order, see \cite{Lav,Lav1}) it follows that in the class $\SScatt$
there are no decreasing sequences of the form $\BL_0 \succ \BL_1 \succ \BL_2 \succ \dots$.
Assuming that $\la \SScatt /\rightleftarrows ,\vartriangleleft \ra\not\in \wfr$ we would have a nonempty set $X\subset \SScatt/\rightleftarrows$,
such that for each $y\in X$ there is $z\in X$ satisfying $z\vartriangleleft y$ and, hence, there would be a decreasing sequence
$y_0\vartriangleright y_1\vartriangleright y_2\vartriangleright \dots$ . Choosing $\BL _i\in y_i$, for $i\in \o$, we would obtain a sequence
$\BL_0 \succ \BL_1 \succ \BL_2 \succ \dots$, which is impossible. So, $\la \SScatt /\rightleftarrows ,\vartriangleleft \ra\in \wfr$.

Let $\theta :\CX \rightarrow \SScatt /\rightleftarrows$
be given by $\theta (\X _i)=[\X _i]_\rightleftarrows$. If $\X _i \preccurlyeq _m \X _j$, then, since monomorphisms of linear orders are embeddings,
$\X _i \hookrightarrow \X _j$, and, hence,  $[\X _i]\trianglelefteq [\X _j]$.
Thus (\ref{EQB044}) is true, $\theta\in \M (\CX )$, and we apply Theorem \ref{TB050}.
($I^\theta _{[\BL ]_\rightleftarrows}= \{ i\in I : \X _i \rightleftarrows \BL \}$, $\BL\in \SScatt$, are finite sets.)
\kdok
\begin{ex}\rm
If $\X _i$, $i\in I$, are arbitrary linear orders of size $\leq \o$ such that the sequence $\la [\X _i]_\rightleftarrows :i\in I\ra$ is finite-to-one,
then $\bigcup _{i\in I}\X_i$ is a reversible poset.
This follows from Proposition \ref{TB053}, since $\BL\in \SScatt$, for each countable linear order $\BL$.
We note that, in that case,
 $\X _i\in \Scatt$, for all except finitely many $i\in I$
(because each countable non-scattered linear order is bi-embedable with $\Q$).
\end{ex}
\paragraph{CSB chains of a limit type}
We recall that the {\it order type} of a linear order $\BL$ is the class $\otp (\BL )=[\BL ]_{\cong}$ of all linear orders isomorphic to $\BL$.
$\BL$ will be called {\it Cantor-Schr\"{o}der-Bernstein (for embeddings)}  iff for each linear order $\BL'$
satisfying $\BL'\rightleftarrows \BL$ we have $\BL'\cong \BL$, that is, $[\BL ]_{\rightleftarrows }=[\BL ]_{\cong}$.
So, by Proposition \ref{TB053} we have
\begin{cor}\label{TB055}
{\rm ($\ast$)} If $\X_i$,  $i\in I$, are $\s$-scattered CSB linear orders, then\\[-2mm]

 $\la \otp (\X _i) :i\in I\ra$ is a finite-to-one sequence\footnote{that is,
there is no infinite $J\subset I$ such that $\X _i \cong \X _j$, for all $i,j\in J$.}
 $\Rightarrow$ $\bigcup _{i\in I}\X_i$ is a reversible poset.
\end{cor}
A scattered linear order $\BL$ is said to be {\it of a limit type}  iff $\BL \rightleftarrows \sum _{s\in S} \BL _s$, where
$\S\in \Scatt$ and
$\BL _s\cong \o$ or $\BL _s\cong\o^*$, for each $s\in S$.\footnote{$\BL\in \Scatt$ is of a limit type
iff $\BL$ has no points which are left fixed under every $f\in \Emb (\BL)$,
iff no ${\mathbf 1}$ appears in the expression of $\BL$ as a minimal sum of hereditarily additively indecomposable linear orders; see \cite{Lav1}, p.\ 112.}
For example, the linear order $\BL= \o\o^* +1$ is of a limit type because $\BL\rightleftarrows \o\o^*$ but, since $\BL\not\cong \o\o^*$, it is not
CSB. Successor ordinals are CSB,  but not of a limit type; limit ordinals are CSB of a limit type.

In order to describe CSB chains of a limit type let
$\CW$ denote the class of well orders,
$\L $ the class of well orders isomorphic to limit ordinals,
$\CZ$ the class of linear orders isomorphic to $\o ^\theta \o ^* +\o ^{\d }$, where $\theta$ and $\d$ are ordinals satisfying $1\leq\theta<\d$.
$\CW ^*$, $\L ^*$, and $\CZ ^*$ will denote the classes of the inverses of elements of $\CW $, $\L $, and $\CZ $, respectively.
(Clearly, $(\o ^\theta \o ^* +\o ^{\d })^*=(\o ^\d)^* + (\o ^\theta)^* \o$.)
\begin{te}\label{T8135}
(a) A linear order is CSB of a limit type iff it is isomorphic to a finite sum of linear orders from
$\L \cup \L ^* \cup \CZ \cup \CZ ^*$.

(b) If \,{\rm ($\ast$)} holds and if $\,\X_i$, $i\in I$, are CSB linear orders of a limit type, then\\[-3mm]

$\bigcup _{i\in I}\X_i$ is a reversible poset
$\;\Leftrightarrow\;$ $\la \otp (\X _i) :i\in I\ra$ is a finite-to-one sequence.
\end{te}
\dok
(a) From recent results of Laflamme, Pouzet, and Woodrow, (see \cite{Laf}) it follows that a scattered linear order is CSB iff it is isomorphic to a finite sum of linear orders from
$\CW \cup \CW ^* \cup \CZ \cup \CZ ^*$.

Let $\BL$ be CSB linear order of a limit type, presented as a sum $\BL =\BL _1 + \dots +\BL _n$,
where $\BL _i \in \CW \cup \CW ^* \cup \CZ \cup \CZ ^*$, for $i\leq n$. Let, in addition, this is a presentation of $\BL$
with the minimal number of summands from $\CW \cup \CW ^* \cup \CZ \cup \CZ ^*$
(for example, $\o +\o \in \CW$, but this linear order can be presented as a sum of finitely many
elements of $\CW$ in infinitely many ways).

Suppose that $i\leq n$ and $\BL _i \in \CW \setminus \L$, that is $\BL _i =\A+\B$,
where $\A \cong \g\in \L \cup \{ 0\}$, $B=\{ b_0, \dots, b_{k-1}\}$, for some $k\in \N$, and $b_0<b_1< \dots <b_{k-1}$.
Since $\BL$ is of a limit type,
$b_0$ belongs to a convex part $C$ of $\BL$ such that  $C\cong\o$ or $C\cong\o^*$.

If $C\cong\o$, then $C_1:=(b_{k-1}, \infty )_\BL \cap C$ is a convex part of $\BL$ of type $\o$
and $b_k:= \min C_1 = \min \BL _{i+1}$, which implies that
$\BL _{i+1}\in \CW $, because the linear orders from $(\CZ \cup \CZ ^* \cup \CW ^*)\setminus \CW$ have no minimum.
Thus $\BL _i +\BL _{i+1}\in \CW $,
which contradicts the minimality of $n$.

If $C\cong\o^*$, then $\g=0$, $i>1$, $\BL _{i-1} \not\in \CZ \cup \CZ ^*$ (since the linear orders from $\CZ \cup \CZ ^*$ do not have a largest element)
 and, hence, $\BL_{i-1} \in \CW ^*$ and $\BL_{i-1} + \BL_i \in \CW ^*$, which contradicts the minimality of $n$ again.

So, $\BL _i \in \CW $ implies that $\BL _i \in  \L$ and, similarly, $\BL _i \in \CW ^*$ implies that $\BL _i \in  \L^*$.

Conversely, let $\BL =\BL _1 + \dots +\BL _n$,
where $\BL _i \in \L \cup \L ^* \cup \CZ \cup \CZ ^*$, for $i\leq n$, and $n$ is the minimal number of summands.
It is well known that the binary relation $\sim$ on $L$ defined by: $x\sim y $ iff $|[\min\{x,y\},\max\{x,y\} ]|<\o$, is an equivalence relation
and (see \cite{Rosen}, p.\ 71) $\BL =\sum _{t\in T}\BL _t$, where $\T\in \Scatt$, $\{ L_t :t\in T\}=L/\sim $ and
$\otp (\BL _t)\in \N \cup \{ \o ,\o ^* ,\zeta\}$, for each $t\in T$ (where $\zeta :=\otp (\Z )$). So, since
each $x\in L$ belongs to a convex subset $C$ of $L$, which is either isomorphic to $\o$ (if $x\in L_i$ and $\BL_i \in \L \cup\CZ$)
or to $\o ^*$ (if $x\in L_i$ and $\BL_i \in \L ^* \cup\CZ ^*$),
we have $\otp (\BL _t )\not\in \N$, for all $t\in T$. Thus $\BL$ can be presented as a sum of linear orders isomorphic to
$\o $, $\o ^*$, or $\zeta =\o ^* +\o $ and, hence,  it is of a limit type.

(b) The implication ``$\Leftarrow$" follows from  Corollary \ref{TB055}.

Suppose that the sequence $\la \otp (\X _i) :i\in I\ra$ is not finite-to-one
and $J=\{ i_k :k\in \o\}\subset I$, where $i_k \neq i_l$ and $\X _{i_k} \cong \X _{i_l}$, for different $i,j\in J$.
By the assumptions we have $\X _{i_0}\cong \sum _{s\in S} \BL _s$, where $\BL _s\cong \o$ or $\BL _s\cong\o^*$.
For $s\in S$ let $L_s =\{ a^s_0, a^s_1, a^s_2, \dots  \}$ be an enumeration such that $a^s_0 < a^s_1 < a^s_2< \dots$, if $\BL _s\cong \o$,
and $a^s_0 > a^s_1 > a^s_2> \dots$, if $\BL _s\cong \o^*$.
Then defining $A_0:=\{ a^s_{2n}: s\in S \land n\in \o\}$ and $A_1:=\{ a^s_{2n+1}: s\in S \land n\in \o\}$ we have
$\A_0\cong \A_1 \cong \X _{i_0}$ and $\{ A_0, A_1\}$ is a partition of $X _{i_0}$. Let $f\in \Sur (I)\setminus \Sym (I)$ be defined by
$f(i_0)=f(i_1)=i_0$, $f(i_k)=i_{k-1}$, for $k\in \N$, and $f(i)=i$, for $i\in I\setminus \{ i_k :k\in \o\}$.
Then $X _{i_0}$ is partitioned into copies of $X _{i_0}$ and $X _{i_1}$ and, by Corollary \ref{TB052}, the poset $\bigcup _{i\in I}\X_i$ is not reversible.
\kdok
We note that, in particular, Theorem \ref{T8135}(b) gives a characterization of reversibility
in the class of posets of the form $\X=\bigcup _{i\in I}\X_i$, where $\X_i\in \L \cup \L ^*$, for $i\in I$.
The corresponding characterization, when $\X_i\in \CW \cup \CW ^*$, for $i\in I$, is given in \cite{KuMoWO}.
\paragraph{Acknowledgments}
This research was supported by the Ministry of Education and Science of the Republic of Serbia (Project 174006).

\end{document}